\documentclass[12pt]{article}
\usepackage{amssymb}
\def\part#1{\frac{\partial\phantom{q}}{\partial#1}}
\newtheorem{thm}{Theorem}
\newtheorem{prp}[thm]{Proposition}

\newenvironment{proof}[1][Proof]{\par\noindent{\bfseries\itshape #1:}\enspace}
  {\medbreak\medskip}
\newcommand{\lie}[1]{\mathfrak{#1}}
\def\End{\mathop{\rm End}\nolimits}

\def\im{\mathop{\rm Im}\nolimits}

\def\tr{\mathop{\rm tr}\nolimits}

\newcommand{\R}{\mathbf{R}}
\newcommand{\C}{\mathbf{C}}
\newcommand{\Z}{\mathbf{Z}}

\begin{document}
\title{$L^2$-cohomology of  hyperk\"ahler quotients}
\author{Nigel Hitchin\\[5pt]
\itshape  Mathematical Institute\\
\itshape 24-29 St Giles\\
\itshape Oxford OX1 3LB\\
\itshape England\\
 hitchin@maths.ox.ac.uk}
\maketitle
\section{Introduction}
There has been some recent interest in the study of $L^2$ harmonic forms on certain non-compact moduli spaces occurring in gauge theories. In this note, we use a result of Jost and Zuo \cite{JZ} to prove some of the properties physicists expect of these forms.

 Jost and Zuo's theorem (adapting an earlier idea of Gromov \cite{Grom}) states that if the K\"ahler form $\omega$ on a complete K\"ahler manifold  satisfies $\omega =d\beta$, where $\beta$ is a one-form of linear growth, then the only $L^2$ harmonic forms lie in the middle dimension. An application of the same  argument shows further that if $G$ is a connected Lie group of isometries on a complete Riemannian manifold    generated by Killing vector fields of linear growth, then $G$ acts trivially on the space of $L^2$ harmonic forms.

Since $2i\omega=\sum dz_j\wedge d\bar z_j=d(\sum z_j\wedge d\bar z_j)=d\beta$,  Euclidean space ${\C}^n$ is the model for the growth conditions required. Moreover, Killing vectors on ${\C}^n$ are also of linear growth. Our applications are principally to hyperk\"ahler quotients of flat Euclidean spaces, where these properties are inherited. We extend the finite-dimensional arguments  to some infinite-dimensional quotient constructions, notably moduli spaces of Higgs bundles on Riemann surfaces and monopoles on ${\R}^3$. For the latter we prove some of the S-duality predictions of Sen \cite{Sen}.

The author wishes to thank J. Jost and J. H. Rawnsley for useful comments.

\section{Linear growth}

First recall how Hodge theory works on a complete non-compact Riemannian manifold \cite{dR}.  If $\Omega^p_{(2)}$ denotes the Hilbert space of $L^2$ $p$-forms, then the $L^2$-cohomology group $\bar H^p_{(2)}$ is defined as the quotient of the space of closed ${ L}^2$ $p$-forms by the {\it closure} of the space
$$d\Omega^{p-1}_{(2)}\cap \Omega^p_{(2)}$$
It is a theorem that on a complete manifold any harmonic $ L^2$ $p$-form is closed and coclosed and so represents a class in $\bar H^p_{(2)}$. The Hodge decomposition theorem then implies that there is a unique such representative.
\vskip .25cm
We first reproduce  for the reader's benefit the proof  of the theorem of Jost and Zuo:
\begin{thm} {\cite{JZ}} Let $M$ be a complete oriented Riemannian manifold and let $\alpha=d\beta$ be a  $p$-form such that
$$\Vert \alpha(x) \Vert \le c,\qquad \Vert \beta(x) \Vert \le c' \rho(x_0,x)+c''$$
where $\rho(x_0,x)$ is the Riemannian distance from a  point $x_0\in M$ and $c,c',c''$ are constants. Then for each $ L^2$-cohomology class $[\eta]\in \bar H^q_{(2)}(M)$,
$$[\alpha \wedge \eta]=0\in \bar H^{p+q}_{(2)}(M)$$
\end{thm}
Note that the inequality on $\beta$ is what we mean by ``linear growth''.

\begin{proof} Let $B_r$ be the ball in $M$ with centre $x_0$ and radius $r$. Take a smooth function $\chi_r:M\rightarrow {\R}^{+}$ with
\begin{eqnarray*}
0\le \chi_r(x)&\le& 1\\
\chi_r(x)&=&1\,\,{\rm for}\,\,x\in B_r\\
\chi_r(x)&=&0\,\,{\rm for}\,\,x\in M\backslash B_{2r}\\
\Vert d\chi_r(x) \Vert& \le &K/\rho(x_0,x)\,\,{\rm for}\,\,x\in B_{2r}\backslash B_r
\end{eqnarray*}
Such a function may be obtained by smoothing the function $f(\rho(x_0,x))$ where $f(\rho)=1$ for $\rho \le r$, $f(\rho)=2-\rho/r$ for $r\le \rho\le 2r$ and $f(\rho)=0$ for $\rho\ge 2r$.

The form $d(\chi_r\beta \wedge \eta)$ has compact support, so $d(\chi_r\beta \wedge \eta)\in \Omega^{p+q}_{(2)}$. We want to show that as $r\rightarrow \infty$ these forms converge in ${ L}^2$ to $\alpha\wedge \eta$. Consider
\begin{equation}
d(\chi_r\beta \wedge \eta)=d\chi_r\wedge\beta \wedge \eta +\chi_r\alpha \wedge \eta
\label{terms}
\end{equation}
Since $\Vert \alpha(x) \Vert \le c $, and $\eta \in L^2$, then $\alpha \wedge \eta \in L^2$ and hence
\begin{equation}
\int_M\Vert \alpha \wedge \eta\Vert^2 =\lim_{r\rightarrow \infty}\int_{B_r}\Vert \alpha \wedge \eta \Vert^2
\label{lim1}
\end{equation}
As $\chi_r$ vanishes outside $B_{2r}$ and is identically $1$ on $B_r$, we have 
\begin{equation}
\lim_{r\rightarrow \infty}\int_M\vert \chi_r\vert^2\Vert \alpha \wedge \eta\Vert^2=\lim_{r\rightarrow \infty}\int_{B_r}\Vert \alpha \wedge \eta\Vert^2 + \lim_{r\rightarrow \infty}\int_{B_{2r}\backslash B_r}\vert \chi_r\vert^2\Vert \alpha \wedge \eta\Vert^2
\label{lim2}
\end{equation}
But
$$\int_{B_{2r}\backslash B_r}\vert \chi_r\vert^2\Vert \alpha \wedge \eta\Vert^2\le \int_{B_{2r}\backslash B_r}\Vert \alpha \wedge \eta\Vert^2$$
and since $\alpha \wedge \eta \in { L}^2$, the right hand side tends to zero as $r\rightarrow \infty$, thus so does the left hand side. From (\ref{lim1}) and (\ref{lim2}) we see that $\chi_r\alpha\wedge \eta$  converges in ${ L}^2$ to $\alpha \wedge \eta$.
\vskip .25cm
Now $d\chi_r$ vanishes on $B_r$ and outside $B_{2r}$, and on the annulus in between we have the  estimates  $\Vert d\chi_r(x) \Vert\le K/\rho(x_0,x)$ and  $\Vert \beta(x) \Vert \le c' \rho(x_0,x)+c''$. Thus 
$$\int_M\Vert d\chi_r\wedge \beta\wedge \eta\Vert^2\le {\mathrm const.}\int_{B_{2r}\backslash B_r}\Vert \eta \Vert^2$$
This again converges to zero as $r\rightarrow \infty$ since $\eta \in { L}^2$. We thus have  convergence of both terms on the right hand side of  (\ref{terms}) and consequently $d(\chi_r\beta\wedge \eta)$ converges in ${ L}^2$ to $\alpha \wedge \eta$. Hence $\alpha \wedge \eta$ lies in the closure of $d\Omega^{p-1}_{(2)}\cap \Omega^p_{(2)}$ and its $L^2$-cohomology class vanishes.
\end{proof}
\vskip .25cm
\begin{thm} Let $M$ be a complete K\"ahler manifold of complex dimension $n$ such that the K\"ahler form satisfies $\omega=d\beta$ where 
$$\Vert \beta(x)\Vert  \le c' \rho(x_0,x)+c''$$
then  all ${ L}^2$ harmonic $p$-forms for $p\ne n$ vanish.
\end{thm}
\begin{proof} Since $\omega$ is covariant constant, $\Vert \omega \Vert$ is constant. Thus
 from Theorem 1, the linear growth of $\beta$ implies that the map $L:H^p_{(2)}(M)\rightarrow H^{p+2}_{(2)}(M)$ defined by
$L([\eta])=[\omega \wedge \eta]$ is zero. By Hodge theory this means that if $\eta$ is an $L^2$   harmonic $p$-form, then the ${ L}^2$ harmonic $(2n-p)$-form $\omega^{n-p}\wedge \eta$ vanishes for $p<n$. But by linear algebra, the map $\eta \mapsto \omega^{n-p}\wedge \eta$ is an isomorphism, hence the only non-zero harmonic forms occur when $p=n$.
\end{proof}

Using the  technique of Theorem 1 we now prove the following:
\begin{thm}  Let $M$ be a complete oriented Riemannian manifold and let $G$ be a connected  Lie group of isometries such that the Killing vector fields $X$ it defines  satisfy
$$\Vert X(x)\Vert \le c' \rho(x_0,x)+c''$$
Then  each ${L}^2$-cohomology class is fixed by $G$.
\end{thm}
\begin{proof} 
The group $G$ acts unitarily on the Hilbert space of ${ L}^2$ harmonic forms. This may possibly be infinite-dimensional. Nevertheless, one knows that for a unitary representation the space of analytic  vectors -- the ones for which $g\mapsto g\cdot v$ is analytic -- is dense \cite{N}. If $\eta$ is an $L^2$ harmonic form which lies in this  subspace, then it has a well-defined $L^2$ Lie derivative
$${\cal L}_X\eta=d(\iota(X)\eta)+\iota(X)d\eta=d(\iota(X)\eta)$$
which is also a harmonic form. As in Theorem 1, we write
\begin{equation}
d(\chi_r\iota(X)\eta)=d\chi_r\wedge\iota(X)\eta +\chi_rd(\iota(X)\eta)
\label{terms1}
\end{equation}
We now proceed as before: $\chi_r d(\iota(X)\eta)$ converges in ${ L}^2$ to $d(\iota(X)\eta)$, and from the linear growth of $X$, $d\chi_r\wedge\iota(X)\eta$ converges in ${ L}^2$ to zero. Thus from (\ref{terms1}), $d(\chi_r\iota(X)\eta)$ converges to $d(\iota(X)\eta)$, whose $L^2$-cohomology class is therefore zero. Hence $G$ acts trivially on a dense subspace of ${ L}^2$ harmonic forms and by continuity is trivial on the whole space.
\end{proof}

\section {Hyperk\"ahler quotients}
We shall  apply the above results to certain hyperk\"ahler manifolds. Recall that a hyperk\"ahler manifold is a Riemannian manifold with compatible covariant constant complex structures $I,J,K$ satisfying the quaternionic identities
$$I^2=J^2=K^2=IJK=-1$$
The corresponding K\"ahler forms are $\omega_1,\omega_2,\omega_3$. If we fix the complex structure $I$, then $\omega^c=\omega_2+i\omega_3$ is a holomorphic symplectic form.
\vskip .25cm
On a hyperk\"ahler manifold the Ricci tensor vanishes. This already means that Bochner-type vanishing theorems for some $L^2$ harmonic forms can be applied, using the modification of this approach due to Dodziuk \cite{Dod}. In particular, there are no $L^2$ harmonic $1$-forms or $(p,0)$ forms. If one of the K\"ahler forms satisfies the conditions of Theorem 1, then we have more:

\begin{thm} Let $M$ be a complete hyperk\"ahler manifold of real dimension $4k$ such that one of the  K\"ahler forms  $\omega_i=d\beta$ where $\beta$ has linear growth. Then any $L^2$ harmonic form is primitive and of type $(k,k)$ with respect to all complex structures.
\end{thm}
\begin{proof} The three $2$-forms $\omega_1,\omega_2,\omega_3$ define by exterior multiplication  three commuting operators $L_1,L_2,L_3$ on the algebra of differential forms. There are three  adjoints $\Lambda_1,\Lambda_2,\Lambda_3$, and it is a matter of linear algebra to show that these satisfy the following  commutation relations:
\begin{equation}
[L_1,\Lambda_2] =  [\Lambda_1,L_2] =  -\sigma_3  \label{so5} 
\end{equation}
and similar ones by cyclic permutation.
Here $\{\sigma_1,\sigma_2,\sigma_3\}$ is the standard basis of the Lie algebra $\mathfrak{su}(2)$ satisfying $[\sigma_1,\sigma_2]=2\sigma_3$ etc. Its action is induced from the action of the unit quaternions on the exterior algebra.  
(As noted by  Verbitsky \cite{Verb}, $L_i$ and $\Lambda_j$  generate an action of the Lie algebra $\lie{sp}(1,1)\cong \lie{so}(5,1)$).

From Theorem 2, each $L^2$ harmonic form $\eta$  lies in the middle dimension $2k$. But  $L_i$ and $\Lambda_i$ commute with the Laplacian, and map to $2k+2$ and $2k-2$ forms respectively, so  we must have for each $i=1,2,3$
$$L_i\eta=0=\Lambda_i\eta.$$
But then  (\ref{so5}) implies that $\sigma_i \eta=0$. For a complex structure $I$ corresponding to $\sigma_i$, the eigenspaces of $\sigma_i$ are the forms of type $(p,q)$ with eigenvalue $i(p-q)$, so if $\sigma_i \eta=0$, then $p=q=k$. The condition $\Lambda_i\eta=0$ is the statement that the form is primitive, so we see that $\eta$ is of type $(k,k)$ and primitive with respect to all complex structures. In Riemannian terms this implies that $\eta$ is anti-self-dual if $k$ is odd and self-dual if $k$ is even.
\end{proof}
\vskip .25cm
Interestingly, hyperk\"ahler manifolds provide a good source of examples where the linear growth conditions hold. We can frequently write one of the K\"ahler forms in a canonical way as $\omega_i=d\beta$ if there is a non-trivial Killing field $X$. Any Killing vector field acts by Lie derivative on the space of covariant constant $2$-forms, preserving the metric. If the hyperk\"ahler manifold is irreducible, this space is spanned by $\omega_1,\omega_2,\omega_3$.  If all are annihilated by $X$, then the action is called triholomorphic. The  alternative  case is relevant for us:  for some orthonormal basis
\begin{equation}
{\cal L}_X\omega_1=0,\quad {\cal L}_X\omega_2=\omega_3,\quad {\cal L}_X\omega_3=-\omega_2
\label{omegas}
\end{equation}
then 
$$\omega_3={\cal L}_X\omega_2=d(\iota(X)\omega_2)+\iota(X)d\omega_2=d(\iota(X)\omega_2)$$
since $d\omega_2=0$. Thus if $X$ acts in this way and has linear growth we can apply Theorem 2 with $\beta=\iota(X)\omega_2$, and also Theorem 3 of course. 

This situation occurs frequently when forming hyperk\"ahler quotients. If $G$ is a Lie group acting on a hyperk\"ahler manifold $M$ preserving the three K\"ahler forms $\omega_1,\omega_2,\omega_3$ then we have three moment maps giving a single function $\mu=(\mu_1,\mu_2,\mu_3)$:
$$\mu: M \rightarrow {\lie g}^*\otimes {\R}^3$$
Suppose $\mu^{-1}(0)$ is smooth then the induced metric is $G$-invariant and descends to the quotient. The hyperk\"ahler quotient construction \cite{HKLR} is the observation that this quotient metric is again hyperk\"ahler.

Suppose we take $A$ to be a complex $n$-dimensional affine space with a flat Hermitian metric and $G$ to be a Lie group of unitary isometries. The cotangent bundle
$$T^*A\cong A\times {\C}^n$$
 is a flat  hyperk\"ahler manifold, with $\omega^c$ the canonical symplectic form. The natural action of the group $G$ preserves $\omega^c$ and also the  hermitian form induced from that of $A$ and so preserves the three K\"ahler forms. On the other hand there is a ${\C}^*$ action given by scalar multiplication of the cotangent vectors by a non-zero complex number $\lambda$. Since the canonical form transforms as $\omega^c\mapsto \lambda\omega^c$, the action of the circle  $S^1\subset {\C}^*$ defines a vector field $X$ with the properties (\ref{omegas}). This commutes moreover with the action of $G$. 

For a vector field $Y$ on a manifold, there is a natural choice for the moment map of the canonical lift to the cotangent bundle: $\mu_Y(\theta_x)=\theta_x(Y_x)$. Make this choice for $\mu^c=\mu_2+i\mu_3$, and an arbitrary choice for $\mu_1$. Then the circle acts on $\mu^{-1}(0)$, commutes with $G$ and hence descends to an action satisfying (\ref{omegas}) on the hyperk\"ahler quotient. 
\vskip .25cm
Most quotients produced this way are complete for  $\mu^{-1}(0)$ is a closed submanifold of a complete manifold and hence complete. Certainly if $G$ is compact the quotient is then complete. We can also obtain growth estimates on Killing vector fields:
\begin{prp} Let $H$ be a group of isometries of $T^*A$ which preserves $\mu^{-1}(0)$ and normalizes $G$. Then $H$ acts isometrically on the quotient and the corresponding Killing vector fields have linear growth.
\end{prp}
\begin{proof}
 Let $Y$ be a Killing vector field on $T^*A$ generated by the action of $H$. Then $Y$ is tangential to $\mu^{-1}(0)$. Let  $\bar Y$ be the corresponding vector field on the quotient $\bar M$, then the metric on $T_{\bar x}\bar M$ is the induced inner product on the horizontal space in $T_x(\mu^{-1}(0))$: the orthogonal complement of the tangent space of the orbit of $G$. Thus at $\bar x \in \bar M$
\begin{equation}
(\bar Y,\bar Y)_{\bar x}= (Y_H,Y_H)_x\le (Y,Y)_x
\label{less}
\end{equation}
The distance $\rho(\bar x,\bar x_0)$ between points $\bar x,\bar x_0\in\bar M$ is the length in $\mu^{-1}(0)$ of the horizontal lift from $x_0$ to $x$ of a geodesic in $\bar M$, and this is greater than or equal to the straight line distance between $x$ and $x_0$, thus
\begin{equation}
 \rho(\bar x,\bar x_0)\ge \Vert x-x_0\Vert
 \label{dist}
 \end{equation}
But the vector field $Y$ is defined by a group of  affine transformations of a flat space and so is of the form
$$(\sum A_{ij}x_i+b_j)\frac{\partial}{\partial x_j}$$
which has linear growth, so from (\ref{less}), 
$$\Vert \bar Y \Vert_{\bar x}\le c'\Vert x-x_0\Vert +c''$$
and  using (\ref{dist}),
$$\Vert\bar Y \Vert_{\bar x}\le c'\rho(\bar x,\bar x_0) +c''$$
\end{proof}
\vskip .25cm
There are many examples in the literature of such complete finite-dimensional quotients, where we can address the question of existence of $L^2$ harmonic  forms.

\subsection{The Taub-NUT metric}

In this case we take $G={\R}$ acting on $A=\C^2$ by
$$(z_1,z_2)\mapsto (e^{it}z_1,z_2+t)$$
Concretely, we have coordinates $z_1,z_2,w_1,w_2$ on $T^*A$, and the action is
$$(z_1, z_2, w_1, w_2)\mapsto (e^{it}z_1,z_2+t, e^{-it}w_1,w_2)$$
The complex moment map is 
$$\mu^c=iz_1w_1+w_2$$
and the real moment map
$$\mu_1=\vert z_1\vert^2-\vert w_1 \vert^2+\im z_2$$
In each ${\R}$-orbit there is a unique point with $z_2$ imaginary. The moment map equations $\mu_1=\mu^c=0$ then define $w_2$ and $z_2$ in terms of $z_1,w_1$. The hyperk\"ahler quotient is  ${\C}^2$, with coordinates $z_1,w_1$. Its metric is   the complete Taub-NUT metric.

 Theorem 4 now tells us that the only $L^2$ harmonic forms on Taub-NUT space are of type $(1,1)$ and  primitive for all complex structures $I,J,K$. In four dimensions ($k=1$ in Theorem 4) this means that they are anti-self-dual. 
\vskip .25cm
To see an $L^2$ harmonic 2-form on this space, we can follow \cite{Gib}. The circle action
$$(e^{i\theta}z_1, z_2, e^{-i\theta}w_1, w_2)$$
commutes with $G$ and preserves all K\"ahler forms, and so descends to the Taub-NUT space with the same property. This gives the Killing vector field $X=\partial /\partial \tau$ which expresses the Taub-NUT metric in its more usual form
$$g=V(dx_1^2+dx_2^2+dx_3^2)+V^{-1}(d\tau +\alpha)^2$$
with $\alpha \in \Omega^1(\R^3)$ and
$$V=1+\frac{m}{r},\qquad d\alpha=\ast dV$$
Here the 1-form dual to $X$ using the metric is
$$\theta=V^{-1}(d\tau +\alpha)$$
and 
$$d\theta=-{V^{-2}}dV\wedge (d\tau +\alpha)+V^{-1}\ast dV$$
which from the form of the metric is anti-self-dual and closed, hence harmonic. Now
$$d\theta\wedge \ast d\theta=-(d\theta)^2=2\frac{1}{V^3}dV\wedge \ast dV\wedge d\tau\simeq \frac{2m^2}{r^4}dx_1\wedge dx_2\wedge dx_3\wedge d\tau$$
and so $d\theta$ is  in $L^2$. It defines a nontrivial $L^2$-cohomology class but as an ordinary cohomology class it is of course trivial.
\vskip .25cm
We shall use Theorem 3 later to give a rigorous proof that this is (up to a constant) the unique $L^2$ harmonic form. There are deformations of products of Taub-NUT spaces which occur as monopole moduli spaces \cite{Gib},\cite{LWY} to which our theorems apply.

\subsection {The Calabi metrics}
Here we take $G=S^1$ acting on $A={\C}^n$ by
$$e^{it}(z_1,\dots,z_n)=(e^{it}z_1,\dots,e^{it}z_n)$$
The induced action on $T^*\C^n\cong \C^{2n}$ is 
$$e^{it}(z_1,\dots,z_n,w_1,\dots,w_n)=(e^{it}z_1,\dots,e^{it}z_n,e^{-it}w_1,\dots, e^{-it}w_n)$$
and this has complex moment map
$$\mu^c=i\sum z_jw_j$$
and real moment map
$$\mu_1=\sum z_j\bar z_j$$
Since $S^1$ is abelian, we can change $\mu_1$ by adding a constant so one choice of  moment map equations $\mu_1=\mu^c=0$ is:
$$\sum z_jw_j=0,\qquad \sum z_j\bar z_j=1$$
The quotient of this is the Calabi metric \cite{Cal} on the cotangent bundle $T^*\C P^{n-1}$.
\vskip .25cm
As argued  by Segal and Selby in \cite{SS}, if on a complete Riemannian manifold, a compactly supported cohomology class defines a  non-trivial  ordinary cohomology class, then there is a non-zero  
 $L^2$ harmonic form representing it. Now the zero section of $T^*\C P^{n-1}$ is Poincar\'e dual to a compactly supported class in the middle dimension. If we evaluate this on the homology class of the zero section we obtain the self-intersection number: the Euler class of the normal bundle. The normal bundle of the zero section of a vector bundle is canonically isomorphic to the bundle itself, which in this case is the  cotangent bundle, whose Euler class is the  Euler characteristic of $\C P^{n-1}$, which is $n$, and in particular non-zero. Thus the compactly supported class is non-trivial in ordinary cohomology and so there exists a corresponding $L^2$ harmonic $(2n-2)$-form.
\vskip .25cm
The Calabi metric has isometry  group $U(n)$, so that Theorem 3 tells us that this form -- and indeed any other $L^2$ harmonic form if such exist -- is invariant by this group.

\section{Moduli spaces}
The main interest in $L^2$ harmonic forms has arisen from the consideration of certain gauge-theoretic hyperk\"ahler moduli spaces. Sometimes, as in \cite{Gib}, these have a finite-dimensional hyperk\"ahler quotient description, but more often than not we need to rely on infinite-dimensional approaches to describe the hyperk\"ahler structure. We give two examples here where we again have a circle action and linear growth conditions which make the Jost-Zuo result applicable.

\subsection{Higgs bundle moduli spaces}

Formally speaking, this is a quotient  like the ones we have just described. Take a compact Riemann surface $\Sigma$ and a  $C^{\infty}$ unitary vector bundle $E$ over $\Sigma$, and let ${\cal A}$ be the space of all unitary connections on $E$. This is an affine space with a constant symplectic form on it: a tangent vector to the space of connections is $ \alpha \in \Omega^1(\Sigma,\End E)$ and then
$$\int_{\Sigma} \tr (\alpha\wedge \beta)$$
defines the form. The $(0,1)$ part $d_A''$ of a connection defines a holomorphic structure on $E$, and this gives the space ${\cal A}$ the structure of a complex affine space. The group ${\cal G}$ of unitary gauge transformations acts on ${\cal A}$ preserving both symplectic and complex structures, and hence the corresponding Hermitian form. As before, we can define a hyperk\"ahler quotient from the action of ${\cal G}$ on the cotangent bundle $T^*{\cal A}$. In this case we have formally
$$T^*{\cal A}\cong {\cal A}\times C^{\infty}(\Sigma, \End E \otimes K)$$
The complex moment map for a pair $(A,\Phi)$ is 
$$\mu^c=d_A''\Phi$$
and the real moment map
$$\mu_1=F_A+[\Phi,\Phi^*]$$
where $F_A$ is the curvature of the connection $A$ . The equations $\mu^c=0,\mu_1=0$ are the Higgs bundle equations, and their quotient by the group of gauge transformations is the moduli space of Higgs bundles. The formal aspects of the foregoing can be made precise (see \cite{H1}) by using some gauge-theoretic analytical results so that, under suitable topological conditions, the $L^2$ inner product gives the moduli space the structure of a complete finite-dimensional hyperk\"ahler manifold. 
\vskip .25cm
There is the same cotangent circle action as before
$$\Phi \mapsto e^{i\theta}\Phi$$
and the vector field induced by this action on the moduli space has linear growth in exactly the same way as in the finite-dimensional case, so we deduce that any $L^2$ harmonic form must be in the middle dimension.
\vskip .25cm
In the case when the bundle $E$ has rank 2 and is of odd degree, Hausel \cite{Hau} studied the image of the compactly supported cohomology in the ordinary cohomology and rather surprisingly showed that it is zero. We cannot therefore assert the existence of $L^2$ harmonic forms by topological means here as for the Calabi metrics, though it is possible that some forms may exist just as they do for the Taub-NUT metric.

\subsection{Monopole moduli spaces}

Here to obtain a hyperk\"ahler quotient we consider a connection $\nabla$ in a principal bundle over ${\R}^3$ and a Higgs field $\phi$, a section of the associated bundle of Lie algebras. The Dirac operator
$$i\nabla_1+j\nabla_2+k\nabla_3-\phi$$
can be considered, with appropriate conditions at infinity, as lying in an affine quaternionic space with an $L^2$ inner product. It is acted on by the group ${\cal G}$ of gauge transformations.

The hyperk\"ahler moment maps for this action are
$$\mu_1=F_{23}-\nabla_1\phi,\quad \mu_2=F_{31}-\nabla_2\phi,\quad \mu_3=F_{12}-\nabla_3\phi$$
and setting $\mu=0$ gives the Bogomolny equations
$$F=\ast \nabla\phi$$
The moduli space of solutions is thus formally a hyperk\"ahler quotient.  There is one complex structure for each direction in ${\R}^3$ and the action of the rotation group $SO(3)$ on ${\R}^3$ gives, for each direction, a circle action whose vector field $X$ is of  the type (\ref{omegas}), fixing only one complex structure. 
\vskip .25cm
We shall show that $X$ has linear growth. To do this, however, we should note that we are not in the simple cotangent bundle  situation any longer. The circle acts non-trivially on the group ${\cal G}$ and moreover its action on $(\nabla, \phi)$ is not algebraic but involves Lie derivatives. To remedy this, we adopt the dual approach through Nahm's equations, restricting ourselves to the most familiar situation for $SU(2)$ monopoles of charge $k$. As in \cite{AH}, analytical results show that the moduli space is a  smooth $4k$-manifold with a complete hyperk\"ahler metric. In this case, thanks to a result of Nakajima \cite{Nak}, the same hyperk\"ahler metric can be seen as a hyperk\"ahler quotient of a space of Nahm matrices as follows.
\vskip .25cm
We  consider the space ${\cal A}$
 of operators $d/ds+B_0+iB_1+jB_2+kB_3$ on functions $f:[0,2]\rightarrow {\C}^k$ with  $B_i:(0,2)\rightarrow {\lie u}(k)$
 where at $s=0$,
 $B_0$ is smooth and for $i=1,2,3$ there is a smooth map $\beta_i:[0,2]\rightarrow {\lie u}(k)$ such that near zero 
 $$B_i=\frac{\rho_i}{s}+\beta_i$$
 for a {\it fixed} irreducible representation of the Lie algebra of $SU(2)$ defined by $\rho_i$. At $s=2$ we
 have the same behaviour.
 \vskip .15cm
 Tangent  vectors $(A_0,A_1,A_2,A_3)$ to this space are  smooth at the
 end-points,  and  using  the group ${\cal G}^0_0$ of smooth maps $g:[0,2]\rightarrow U(k)$ for which
 $g(0)=g(1)=1$,
 and  some  analysis,  we  obtain  a  hyperk\"ahler  metric on the space
 ${\cal B}$
 of solutions to the hyperk\"ahler moment map equations 
  \begin{eqnarray}
B_1'+[B_0,B_1]&=&[B_2,B_3]\nonumber\\
B_2'+[B_0,B_2]&=&[B_3,B_1]\label{beqns}\\
B_3'+[B_0,B_3]&=&[B_1,B_2]\nonumber
\end{eqnarray}
modulo the action of the gauge group ${\cal G}^0_0$. 
\vskip .25cm
This gives the metric, but we need to determine the circle action. Having fixed the
residues, this is less easy to describe, because a rotation involves a compensating
gauge transformation outside ${\cal G}^0_0$. The
infinitesimal version of the action is represented by a vector field
$$X=(\psi'+[B_0,\psi],[B_1,\psi],B_3+[B_2,\psi],-B_2+[B_3,\psi])$$
This is a vector field on the space ${\cal B}$ in
the  infinite-dimensional flat space ${\cal A}$:  we are using the linear
structure of the ambient space to write down tangent vectors.
It must be smooth at the end-points, so to keep the residues fixed we need
$$\psi(0)=\psi(2)=-\rho_1.$$
The K\"ahler form $\omega_1$ on the quotient which is invariant by the induced  action   pulls  back to ${\cal B}$ as the restriction of
the constant symplectic form on ${\cal A}$:
$$ \int_0^2[-\tr(A_0\tilde A_1)+\tr(A_1\tilde A_0)+\tr(A_2\tilde
A_3)-\tr(A_3\tilde A_2)] ds$$
This form is degenerate in the ${\cal G}^0_0$-orbit directions, so the 1-form $\iota(X)\omega_1$ pulls back to
\begin{eqnarray*}
\pi^*\iota(X)\omega_1(A)&=&\int_0^2-\tr(A_0[B_1,\psi])+\tr(A_1(\psi'+[B_0,\psi]))+\tr(A_2(
-B_2+[B_3,\psi]))\\
& &-\tr(A_3(B_3+[B_2,\psi])) ds\\
&=&\int_0^2-\tr([A_0,B_1]\psi)+\tr(A_1\psi')-\tr([B_0,A_1]\psi)-\tr(A_2B_2)\\
& &+\tr([A_2,B_3]\psi)-\tr(A_3B_3)+\tr([B_2,A_3]\psi) ds
\end{eqnarray*}
On the other hand $A$ is tangent to ${\cal B}$, so $A$ satisfies the
linearization of the equations (\ref{beqns}). In particular
\begin{equation}
A_1'+[A_0,B_1]+[B_0,A_1]=[A_2,B_3]+[B_2,A_3]
\label{def}
\end{equation}
Substituting in the formula for $\pi^*\iota(X)\omega_1$ then gives
\begin{eqnarray*}
\pi^*\iota(X)\omega_1(A)&=&\int_0^2(\tr(A_1'\psi)+\tr(A_1'\psi))ds
-\int_0^2\tr(A_2B_2+A_3B_3)ds\\
&=&[\tr(A_1\psi_1)]_0^2-\int_0^2\tr(A_2B_2+A_3B_3)ds\\
\end{eqnarray*}
At first sight it looks as if this integral may not be well-defined, since $B_2$ and $B_3$ have poles at $s=0,2$. In fact looking at the three deformation equations like (\ref{def}) and applying some algebra (as in \cite{H2}) it can be seen that $A_i(0)$ is a scalar, so since $B_i$ has residue $\rho_i$ and $\tr \rho_i=0$, we have a smooth integrand. Moreover $\tr(A_1\psi_1)=c\tr \rho_1=0$, at $s=0,2$, so
$$\pi^*\iota(X)\omega_1(A)=-\int_0^2\tr(A_2B_2+A_3B_3)ds$$
It is clear now from the dependence on $B$ in this formula that $\pi^*\iota(X)\omega_1$ has linear growth in the ambient flat $L^2$  metric and hence by the arguments of Proposition 5, $\iota(X)\omega_1$ and hence $X$ itself has linear growth in the metric on the moduli space.
\vskip .25cm
Let us denote by $M_k$ the $4k$-dimensional moduli space of $SU(2)$ monopoles of charge $k$. Does it have any $L^2$ harmonic forms? The answer is no. To see this we note that translation in ${\R}^3$ induces an isometry of $M_k$. From the Nahm point of view translation by $(x_1,x_2,x_3)$ is defined by
$$(B_0,B_1,B_2,B_3)\mapsto (B_0,B_1+ix_1I,B_2+ix_2I,B_3+ix_2I)$$
This preserves the residues and is clearly of bounded growth, even stronger than linear growth, and so from Theorem 3, any $L^2$ harmonic form is invariant under translation. Now a cyclic $k$-fold covering of $M_k$ splits isometrically
$$\tilde M_k\cong \tilde M_k^0 \times S^1\times {\R}^3$$
where $\tilde M_k^0$ is the space of strongly centred monopoles and the induced action of translation is just Euclidean translation in the ${\R}^3$ factor. Clearly a translation-invariant form cannot be square integrable.
\vskip .25cm
There do exist $L^2$ harmonic forms on $\tilde M_k^0$, however. This is the subject of Sen's conjectures concerning S-duality \cite{Sen}. The manifold  $\tilde M_k^0$ has an isometric ${\Z}/k$ action on it. Let ${\cal H}^p_k$ denote the space of square-integrable harmonic $p$-forms on $\tilde M_k^0$ and decompose ${\cal H}^p_k$ into representation spaces ${\cal H}^p_{k,\ell}$ where the generator of ${\Z}/k$ acts as $e^{2\pi i \ell/k}$. Then  Sen conjectured that:
\begin{itemize}
\item
if $k$ and $\ell$ are coprime then ${\cal H}^p_{k,\ell}=0$ except in the middle dimension $p=2k-2$, in which case it is one-dimensional,
\item
if $k$ and $\ell$ have a common factor, ${\cal H}^p_{k,\ell}=0$ for all $p$
\end{itemize}
Theorem 2 already shows that ${\cal H}^p_{k,\ell}=0$ if $p\ne 2k-2$ thus verifying part of the conjecture. Also, uniqueness in the conjecture implies that the form must be of type $(p,p)$ (since conjugation interchanges harmonic $(p,q)$ forms and $(q,p)$-forms). Moreover it must be primitive since $\Lambda \eta$ is harmonic. These features we have already seen to hold from Theorem 4. The uniqueness in Sen's conjecture also implies that $SO(3)$ acts trivially on ${\cal H}^p_k$ since it has no nontrivial 1-dimensional representations, and this too we have seen  from Theorem 3.

 Segal and Selby in \cite{SS} showed that the image of the compactly supported cohomology in the ordinary cohomology is one-dimensional in precisely the locations predicted by the conjecture, so there certainly exists a non-trivial vector in ${\cal H}^{2k-2}_{k,\ell}$ for $k,\ell$ coprime.
\vskip .25cm
What remains is  still  the question of uniqueness, despite the restrictions achieved above. Let us  see how this can at least be shown for $k=2$ using our results.
\vskip .25cm
The 2-monopole metric \cite{AH} can be put in the form (with $\pi \le \rho\le \infty$):
$$g=f^2d\rho^2+a^2\sigma_1^2+b^2\sigma_2^2+c^2\sigma_3^2$$
using left invariant forms $\sigma_i$ on $SO(3)$, satisfying $d\sigma_1=\sigma_2\wedge \sigma_3$ etc. We know from Theorem 3 that any $L^2$ harmonic form $\eta$ is $SO(3)$-invariant and from Theorem 4 that it is of type $(1,1)$ and primitive with respect to $\omega_1,\omega_2,\omega_3$, hence  anti-self-dual. Thus $\eta$ is a linear combination $\eta=c_1\varphi_1+c_2\varphi_2+c_3\varphi_3$ of $\varphi_1,\varphi_2,\varphi_3$ where
\begin{eqnarray*}
\varphi_1&=&F_1(\rho)\left( d\sigma_1-\frac{fa}{bc}d\rho\wedge \sigma_1\right)\\
\varphi_2&=&F_2(\rho)\left( d\sigma_2-\frac{fb}{ca}d\rho\wedge \sigma_2\right)\\
\varphi_3&=&F_3(\rho)\left( d\sigma_3-\frac{fc}{ab}d\rho\wedge \sigma_3\right)
\end{eqnarray*}
The form $\eta$ is closed, and since
$$d\varphi_1=d\left( F_1(\rho)( d\sigma_1-\frac{fa}{bc}d\rho\wedge \sigma_1)\right)=\left(F'-F\frac{fa}{bc}\right)d\rho\wedge \sigma_2\wedge \sigma_3$$
with similar expressions for the other terms, this implies that each term $\varphi_i$ with non-zero coefficient $c_i$  is closed and so, following \cite{Sen} 
$$F_1(\rho)=const.\exp\left(-\int_{c}^{\rho}\frac{fa}{bc}(r)dr\right)$$
with similar expressions for $F_2,F_3$. As $\rho \rightarrow \infty$, it is known that 
$$f\simeq -1,\quad a\simeq \rho,\quad b\simeq \rho,\quad c\simeq -2$$
(Note in passing that this behaviour also confirms that the vector fields defined by the $SO(3)$ action have linear growth.) It follows that 
$$\frac{fa}{bc}\simeq \frac{1}{2},\quad \frac{fb}{ca}\simeq \frac{1}{2},\quad \frac{fc}{ab}\simeq \frac{2}{\rho^2}$$
Now 
$$\varphi_1\wedge \ast \varphi_1=-\varphi_1\wedge \varphi_1=2F_1^2\frac{fa}{bc}d\rho\wedge \sigma_1\wedge\sigma_2\wedge\sigma_3$$
and this is integrable as $\rho\rightarrow \infty$ as are $\varphi_2\wedge \ast \varphi_2$ and $\varphi_3\wedge \ast \varphi_3$.

At the opposite extreme, as $\rho\rightarrow \pi$,
$$f\simeq -1,\quad a\simeq 2\rho-2\pi,\quad b\simeq \pi,\quad c\simeq -\pi$$
and so
$$\frac{fa}{bc}\simeq \frac{2\rho-2\pi}{\pi^2},\quad \frac{fb}{ca}\simeq \frac{1}{2\rho-2\pi},\quad \frac{fc}{ab}\simeq \frac{1}{2\rho-2\pi}$$
Thus $F_1$ tends to a constant, but
$$\varphi_2\wedge \ast \varphi_2, \varphi_3\wedge \ast \varphi_3\simeq const. \frac{1}{(\rho-\pi)^2}$$
It follows that only $\varphi_1$ lies in $L^2$. We therefore have, up to a constant multiple, a unique $L^2$  harmonic form on $\tilde M^0_2$ as predicted.

An index-theoretic proof of this uniqueness, following \cite{SSZ}, can be found as part of \cite{PB}.
\vskip .25cm
We can use the same approach to prove the uniqueness of the $L^2$ harmonic form on Taub-NUT space described in 3.1. (cf \cite{LWY}). The metric is now of the form
$$g=f^2d\rho^2+a^2(\sigma_1^2+\sigma_2^2)+b^2\sigma_3^2$$
where the equality of the coefficients of $\sigma_1^2$ and $\sigma_2^2$ arises from the extra isometric $S^1$ action. By the same arguments as the 2-monopole space we are led to consider forms $\varphi_1,\varphi_2,\varphi_3$, but now by Theorem 3, the form must be invariant by the extra circle action and hence
$$\varphi_3=F_3(\rho)\left( d\sigma_3-\frac{fb}{a^2}d\rho\wedge \sigma_3\right)$$
We have already constructed an $L^2$ harmonic form in 3.1, so this must be it.
\vskip .25cm
This same argument can also be applied to the Eguchi-Hanson metric, which is the Calabi metric on $T^*{\C}P^1$. It can be written in the form
$$g=\frac{1}{1-(a/r)^4}dr^2+r^2(\sigma_1^2+\sigma_2^2+(1-(a/r)^4)\sigma_3^2)$$
and has $U(2)$-symmetry, so $\varphi_3$ must be the only $L^2$ harmonic form -- the one we determined topologically in 3.2.

\end{document}